\documentclass[12pt]{article}
\usepackage{amsfonts,epsf,amsmath,amssymb,latexsym}
\usepackage{pgf,tikz}
\usepackage{bbm}

\newtheorem{theorem}{\bf Theorem}[section]
\newtheorem{corollary}[theorem]{\bf Corollary}
\newtheorem{lemma}[theorem]{\bf Lemma}
\newtheorem{proposition}[theorem]{\bf Proposition}

\newtheorem{definition}[theorem]{\bf Definition}

\newcommand{\proof}{\noindent{\bf Proof.\ }}
\newcommand{\qed}{\hfill $\square$ \bigskip}

\def\cp{\,\square\,}

\begin{document}

\title{Daisy cubes and distance cube polynomial}

\author{
Sandi Klav\v zar $^{a,b,c}$
\and
Michel Mollard $^{d}$
}

\date{\today}

\maketitle

\begin{center}
$^a$ Faculty of Mathematics and Physics, University of Ljubljana, Slovenia\\
{\tt sandi.klavzar@fmf.uni-lj.si}
\medskip

$^b$ Faculty of Natural Sciences and Mathematics, University of Maribor, Slovenia\\
\medskip

$^c$ Institute of Mathematics, Physics and Mechanics, Ljubljana, Slovenia\\
\medskip

$^d$ Institut Fourier, CNRS Universit\'e Grenoble Alpes, France\\
{\tt Michel.Mollard@univ-grenoble-alpes.fr}
\end{center}

\begin{center}
{\large Dedicated to the memory of our friend Michel Deza}
\end{center}

\begin{abstract}
Let $X\subseteq \{0,1\}^n$. Then the daisy cube $Q_n(X)$ is introduced as the subgraph of $Q_n$ induced by the intersection of the intervals $I(x,0^n)$ over all $x\in X$. Daisy cubes are partial cubes that include Fibonacci cubes, Lucas cubes, and bipartite wheels. If $u$ is a vertex of a graph $G$, then the distance cube polynomial $D_{G,u}(x,y)$ is introduced as the bivariate polynomial that counts the number of induced subgraphs isomorphic to $Q_k$ at a given distance from the vertex $u$. It is proved that if $G$ is a daisy cube, then $D_{G,0^n}(x,y)= C_{G}(x+y-1)$, where $C_G(x)$ is the previously investigated cube polynomial of $G$. It is also proved that if $G$ is a daisy cube,  then $D_{G,u}(x,-x)=1 $ holds for every vertex $u$ in $G$. 
\end{abstract}


\noindent
{\bf AMS Subj.\ Class.\ (2010)}: 05C31, 05C75, 05A15

\section{Introduction}
\label{sec:intro}

In this paper we introduce a subclass of hypercubes, the members of which will be called daisy cubes. This new class contains numerous important classes of graphs such as Fibonacci cubes~\cite{klavzar-2013}, Lucas cubes~\cite{munarini-2001, taranenko-2013}, gear graphs (alias bipartite wheels)~\cite{kirlangic-2009}, and of course hypercubes themselves. 

Our main motivation for the introduction of daisy cubes are the recent investigations of 
Sayg{\i} and E\u{g}ecio\u{g}lu~\cite{saygi-2017+, saygi-2017++} in which they introduced the $q$-cube polynomial and studied it on Fibonacci cubes $\Gamma_n$ and Lucas cubes $\Lambda_n$. This is a bivariate counting polynomial that keeps track of the number of subcubes that are are a given distance from the vertex $0^n$ in $\Gamma_n$ (resp.\ $\Lambda_n$). Many of the results obtained in~\cite{saygi-2017+} and ~\cite{saygi-2017++} present a refinement of the investigation~\cite{klavzar-2012} of the cube polynomial of Fibonacci cubes and Lucas cube. The latter polynomial is a counting polynomial of induced cubes in a graph; it was introduced in~\cite{bresar-2003} and further studied in~\cite{bresar-2006}. 

We proceed as follows. In the rest of this section we introduce some concepts and notation needed in this paper. In the next section we formally introduced daisy cubes, give several examples of them, and deduce some of their basic properties. In particular, daisy cubes admit isometric embeddings into hypercubes. In Section~\ref{sec:distance-cube-poly} we introduce the 
earlier investigated cube polynomial $C_G(x)$ of a graph $G$, and the distance cube polynomial $D_{G,u}(x,y)$. In the main result of the section (Theorem~\ref{th:DfromC}) we prove a somehow surprising fact that the bivariate distance cube polynomial of an arbitrary daisy cube can be deduced from the univariate cube polynomial. More precisely, if $G$ is a daisy cube, then $D_{G,0^n}(x,y)= C_{G}(x+y-1)$. Several consequences of this theorem are also developed. In particular, if $G$ is a daisy cube, then the polynomials $D_{G,0^n}$ and $C_{G}$ are completely determined by the counting polynomial of the number of vertices at a given distance from the vertex $0^n$. In the final section we prove that $D_{G,u}(x,-x)=1$ holds for every vertex $u$ of a daisy cube $G$.

Let $B=\{0,1\}$. If $u$ is a word of length $n$ over $B$, that is, $u = (u_1,\ldots, u_n)\in B^n$, then we will briefly write $u$ as $u_1\ldots u_n$. The {\em weight} of $u\in B^n$ is $w(u) = \sum_{i=1}^n u_i$, in other words, $w(u)$ is the number of $1$s in word $u$. We will use the power notation for the concatenation of bits, for instance $0^n = 0\ldots 0\in B^n$. 

The {\em $n$-cube} $Q_n$ has the vertex set $B^n$, vertices $u_1\ldots u_n$ and $v_1\ldots v_n$ being adjacent if $u_i\ne v_i$ for exactly one $i\in [n]$, where $[n] = \{1,\ldots, n\}$. The set of all $n$-cubes is referred to as {\em hypercubes}. A {\em Fibonacci word} of length $n$ is a word $u = u_1\ldots u_n \in B^n$ such that $u_i\cdot u_{i+1}=0$ for $1\in [n-1]$. The {\em Fibonacci cube} $\Gamma_n$, $n\geq 1$, is the subgraph of $Q_n$ induced by the Fibonacci words of length $n$. A Fibonacci word $u_1\ldots u_n$ is a {\em Lucas word} if in addition $u_1\cdot u_n = 0$ holds. The {\em Lucas cube} $\Lambda_n$, $n\geq 1$, is the subgraph of $Q_n$ induced by the Lucas words of length $n$. For convenience we also set $\Gamma_0 = K_1 = \Lambda_0$. 

If $u$ and $v$ are vertices of a graph $G$, the the {\em interval} $I_G(u,v)$ between $u$ and $v$ (in $G$) is the set of vertices lying on shortest $u,v$-path, that is, $I_G(u,v) = \{w:\ d(u,v) = d(u,w) + d(w,v)\}$. We will also write $I(u,v)$ when $G$ will be clear from the context. A subgraph $H$ of a graph $G$ is {\em isometric} if $d_H(u,v) = d_G(u,v)$ holds for $u,v\in V(H)$. Isometric subgraphs of hypercubes are called {\em partial cubes}. For general properties of these graphs we refer to the books~\cite[Chapter 19]{deza-2010} and~\cite{ovchinnikov-2011}. See also~\cite{albenque-2016, marc-2016} for a couple of recent developments on partial cubes and references therein for additional results. If $H$ is a subgraph of a graph $G$ and $u\in V(G)$, then the \emph{distance $d(u,H)$ between $u$ and $H$} is $\min_{v \in H}({d_G(u,v)})$. Finally, if $G=(V(G),E(G))$ is a graph and $X\subseteq V(G)$, then $\langle X\rangle$ denotes the subgraph of $G$ induced by $X$.  

\section{Examples and basic properties of daisy cubes}
\label{sec:basic}

Let $\le$ be a partial order on $B^n$ defined with $u_1\ldots u_n \le v_1\ldots v_n$ if $u_i\le v_i$ holds for $i\in [n]$. For $X \subseteq B^n$ we define the graph  $Q_n(X)$ as the subgraph of $Q_n$ with  
$$Q_n(X) = \left\langle \{u\in B^n:\ u\le x\ {\rm for\ some}\ x\in X \} \right\rangle$$
and say that $Q_n(X)$ is a {\em daisy cube (generated by $X$)}. 

Vertex sets of daisy cubes are in extremal combinatorics known as {\em hereditary} or {\em downwards closed sets}, see~\cite[Section 10.2]{jukna-2011}. For instance, a result of Kleitman from~\cite{kleitman-1966} (cf.\ \cite[Theorem 10.6]{jukna-2011}) reads as follows: If $X, Y\subseteq B^n$ are hereditary sets, then $|V(Q_n(X))\cap V(Q_n(Y))|\ge |V(Q_n(X))|\cdot |V(Q_n(Y))|/2^n$. 

Before giving basic properties of daisy cubes let us list some of their important subclasses. 
\begin{itemize}
\item If $X = \{1^n\}$, then $Q_n(X) = Q_n$. 
\item If $X = \{u_1\ldots u_n:\ u_i\cdot u_{i+1}=0, i\in [n-1]\}$, then $Q_n(X) = \Gamma_n$.  
\item If $X = \{u_1\ldots u_n:\ u_i\cdot u_{i+1}=0, i\in [n-1], {\rm and}\ u_1\cdot u_n=0\}$, then $Q_n(X) = \Lambda_n$.
\item If $X = \{110^{n-2}, 0110^{n-3}, \ldots, 0^{n-2}11, 10^{n-1}1\}$, then $Q_n(X) = BW_n$ the bipartite wheel also known as a gear graph.
\item If $X = \{u:\ w(u) \le n-1\}$, then $Q_n(X) = Q_n^{-}$ the vertex-deleted cube.
\end{itemize}

The above example which gives an equivalent description of Fibonaci cubes $\Gamma_n$ can be rephrased by saying that $X$ contains all words that do not contain the subword $11$. This can be generalized by defining $X_k$, $k\ge 2$, as the set of words that do not contain $1^k$. In this way more daisy cubes are obtained; in~\cite{liu-1994} these graphs were named generalized Fibonacci cubes. Today, the term ``generalized Fibonacci cubes" is used for a much larger class of graphs as introduced in~\cite{ilic-2012}. However, no additional daisy cubes are obtained via this generalization, no matter whether they are partial cubes or not, cf.~\cite{wei-2016} for the latter aspect of generalized Fibonacci cubes. 

Note that if $x,y\in X$ and $y\le x$, then $Q_n(X) = Q_n(X\setminus \{y\})$. More generally, if $\widehat{X}$ is the antichain consisting of the maximal elements of the poset $(X,\le)$, then $Q_n(\widehat{X}) = Q_n(X)$. Hence, for a given set $X\subseteq B^n$ it is enough to consider the antichain $\widehat{X}$; we call the vertices of $Q_n(X)$ from  $\widehat{X}$ the {\em maximal vertices} of $Q_n(X)$. For instance, let $X=\{u\in B^n: w(u) \le k\}$. Then the maximal vertices of $Q_n(X)$ are the vertices $u$ with $w(u) = k$. In particular, the vertex-deleted $n$-cube $Q_n^-$ can then be represented as
$$Q_n^- = Q_n(\{u:\ w(u) = n-1\})\,.$$

\begin{proposition}
\label{lem:partial-cube}
If $X\subseteq B^n$, then $Q_n(X)$ is a partial cube. 
\end{proposition}

\proof
Let $u, v\in V(Q_n(X))$ and suppose that $u$ and $v$ differ in coordinates $I=\{i_1,\ldots, i_k\}$. Let $I_u = \{i_j\in I:\ u_{i_j} = 1\}$ and $I_v = \{i_j\in I:\ v_{i_j} = 1\}$, so that $I_u\cup I_v$ is a partition of $I$. 

Clearly, $d_{Q_n}(u,v) = k \le d_{Q_n(X)}(u,v)$. Construct now a path $P$ as a concatenation of paths $P_1$ and $P_2$, where the path $P_1$ starts in $u$ and contains vertices obtained from $u$ by changing one by one the coordinates $i_j\in I_u$, while the path $P_2$ then continues from there and contains vertices obtained by changing one by one the coordinates $i_j\in I_v$. Since for each vertex $x$ of $P_1$ we have $x\le u$ and for each vertex $x$ of $P_2$ we have $x\le v$, it follows that $P$ is an $u,v$-path in $Q_n(X)$. The length of $P$ is $k$ and consequently $d_{Q_n(X)}(u,v)\le k$. We conclude that $d_{Q_n(X)}(u,v) = d_{Q_n}(u,v)$.  
\qed

The following observation will be important for our later studies. 

\begin{lemma}
\label{lem:interval}
Let $X\subseteq B^n$. Then 
$$Q_n(X) = \left\langle \bigcup_{x\in \widehat{X}} I_{Q_n}(x,0^n)\right\rangle\,.$$
\end{lemma}

\proof
Let $u\in V(Q_n(X))$. We have already observed that $Q_n(\widehat{X}) = Q_n(X)$, hence there exists a vertex $x\in \widehat{X}$ such that $u\le x$. It is then straightforward that $u\in I_{Q_n}(x,0^n)$. Therefore, $V(Q_n(X)) = V(Q_n(\widehat{X})) \subseteq  \cup_{x\in \widehat{X}} I_{Q_n}(x,0^n)$. 

Conversely, let $u\in \bigcup_{x\in \widehat{X}} I_{Q_n}(x,0^n)$. Then there exists a fixed vertex $x\in \widehat{X}$ such that $u\in I_{Q_n}(x,0^n)$. But then $u\le x$ and consequently $u\in V(Q_n(X))$ so that $\cup_{x\in \widehat{X}} I_{Q_n}(x,0^n)\subseteq V(Q_n(X))$. 

We conclude that $V(Q_n(X)) = \cup_{x\in \widehat{X}} I_{Q_n}(x,0^n)$. 
\qed

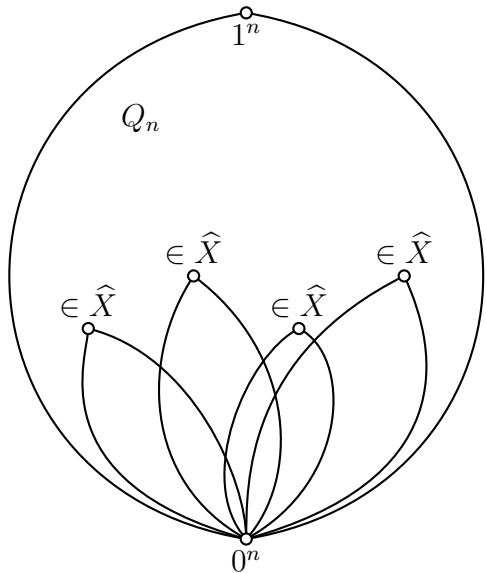
\begin{figure}[ht!]
\begin{center}
\begin{tikzpicture}[scale=0.7,style=thick]
\def\vr{3pt} 
\path (0,0) coordinate (0n);
\path (0,10) coordinate (1n);
\path (-3,4) coordinate (a);
\path (-1,5) coordinate (b);
\path (1,4) coordinate (c);
\path (3,5) coordinate (d);
\draw (0n) .. controls (6,1) and (6,9) .. (1n);
\draw (0n) .. controls (-6,1) and (-6,9) .. (1n);
\draw (0n) .. controls (-4,1.1) and (-3,3.5) .. (a);
\draw (0n) .. controls (-2,1.1) and (-2,3.5) .. (b);
\draw (0n) .. controls (-1,1.1) and (0,3.5) .. (c);
\draw (0n) .. controls (0,1.1) and (0,3.5) .. (d);
\draw (0n) .. controls (0,1.1) and (-1,3.5) .. (a);
\draw (0n) .. controls (1,1.1) and (1,3.5) .. (b);
\draw (0n) .. controls (2,1.1) and (2,3.5) .. (c);
\draw (0n) .. controls (4,1.1) and (3.7,3.5) .. (d);
\draw (0n)  [fill=white] circle (\vr);
\draw (1n)  [fill=white] circle (\vr);
\draw (a)  [fill=white] circle (\vr);
\draw (b)  [fill=white] circle (\vr);
\draw (c)  [fill=white] circle (\vr);
\draw (d)  [fill=white] circle (\vr);
\draw[below] (0n) node {$0^n$};
\draw[below] (1n) node {$1^n$};
\draw (-2,8) node {$Q_n$};
\draw[above] (a) node {$\in \widehat{X}$};
\draw[above] (b) node {$\in \widehat{X}$};
\draw[above] (c) node {$\in \widehat{X}$};
\draw[above] (d) node {$\in \widehat{X}$};
\end{tikzpicture}
	\end{center}
	\caption{A daisy cube}
	\label{fig:daisy}
\end{figure}

Lemma~\ref{lem:interval} is illustrated in Fig.~\ref{fig:daisy}. The figure also gives a clue why the name daisy cubes was selected.

Having Lemma~\ref{lem:interval} in mind we call the vertex $0^n$ the {\em center} of the daisy cube $Q_n(X)$. Let now $X\subseteq B^n$ and let $u\in B^n$ be an arbitrary vertex of $Q_n$. As it is well-known that $Q_n$ is vertex-transitive, there exists an automorphism $\alpha$ of $Q_n$ such that $\alpha(0^n) = u$. Then, in view of Lemma~\ref{lem:interval},  
\begin{equation}
\label{eq:arbitrary-center}
Q_n(X) = \left\langle \bigcup_{x\in \widehat{X}} I_{Q_n}(\alpha(x),u)\right\rangle\,.
\end{equation}
That is, \eqref{eq:arbitrary-center} asserts that for a given daisy cube $Q_n(X)$, an arbitrary vertex of $Q_n$ could an be considered as its center.

\section{Distance cube polynomial}
\label{sec:distance-cube-poly}

For a graph $G$ let $c_k(G)$, $k\ge 0$, be the number of induced subgraphs of $G$ isomorphic to $Q_k$, so that $c_0(G) = |V(G)|$, $c_1(G) = |E(G)|$, and $c_2(G)$ is the number of induced $4$-cycles. The {\em cube polynomial}, $C_G(x)$, of $G$, is the corresponding counting polynomial, that is, the generating function 
$$C_G(x) = \sum_{k\geq 0} c_k(G) x^k\,.$$
Since the cube polynomial is multiplicative on the Cartesian product of graphs and $Q_n$ is the Cartesian product of $n$ copies of $K_2$ for which $C_{K_2}(x)=2+x$ holds, we have
\begin{equation}
\label{eq:cube-poly-for cubes}
C_{Q_n}(x)=(2+x)^n\,.
\end{equation}
In~\cite{saygi-2017+} a $q$-analogue of the cube polynomial of Fibonacci cubes $\Gamma_n$ is considered with the remarkable property that this $q$-analogue counts the number of induced subgraphs isomorphic to $Q_k$ at a given distance from the vertex $0^n$. (For related recent investigations on the number of disjoint hypercubes in Fibonacci cubes see~\cite{gravier-2015, mollard-2017, saygi-2016}.) We now introduce a generalization of this concept to arbitrary graphs as follows.

\begin{definition}
If $u$ is a vertex of a graph $G$, then let $c_{k,d}(G)$, $k,d\ge 0$, be the number of induced subgraphs of $G$ isomorphic to  $Q_k$ at distance $d$ from $u$. The {\em distance cube polynomial} of $G$ with respect to $u$ is
$$D_{G,u}(x,y) = \sum_{k,d\geq 0} c_{k,d}(G) x^ky^d\,.$$
\end{definition}

For the later use we note that 
\begin{equation}
\label{eq:C-from-D}
D_{G,u}(x,1)=C_G(x)\,.
\end{equation} 
We also point out that if $G$ is vertex-transitive, then $D_{G,u}(x,y)$ is independent of $u$.

Recall that the {\em Cartesian product} $G\cp H$ of graphs $G$ and $H$ has the vertex set  $V(G)\times V(H)$ and $E(G\cp H)  = \{(g,h)(g',h'):\ gg'\in E(G) \mbox{ and } h=h', \mbox{ or, } g=g' \mbox{ and }  hh'\in E(H)\}$. If $(g,h)\in V(G\cp H)$, then the {\em $G$-layer}  $G^h$ through the vertex $(g,h)$ is the subgraph of $G\cp H$ induced by the vertices $\{(g',h):\ g'\in V(G)\}$. Similarly, the {\em $H$-layer} $^gH$ through $(g,h)$ is the subgraph of $G\cp H$ induced by the vertices $\{(g,h'):\ h'\in V(H)\}$. The distance cube polynomial is multiplicative on the Cartesian product, just as it is the cube polynomial. 

\begin{proposition}
\label{prp:cart-prod-multiplicative}
If $G$ and $H$ are graphs and $(g,h)\in V(G\cp H)$, then 
$$D_{G\Box H,(g,h)}(x,y) = D_{G,g}(x,y)D_{H,h}(x,y)\,.$$
\end{proposition}

\proof
The result follows from the following fact: If $Q=Q_d$ is a subgraph of $G\cp H$, then for some $k\in [d-1]$ we have $p_G(Q) = Q_k$, $p_H(Q) = Q_{d-k}$, and hence $Q = p_G(Q) \cp p_H(Q)$. (Here $p_G$ and $p_H$ are the projection maps from $G\cp H$ onto the factors $G$ and $H$, respectively).

To prove the above fact consider $(g,h)\in V(Q)$. Let $(g_1,h),\ldots, (g_k,h)$ and $(g,h_1),\ldots, (g,h_{d-k})$ be the neighbors of $(g,h)$ in $Q$. Then by the unique square property of the Cartesian product (see~\cite[Lemma 6.3]{hammack-2011}) we infer that the vertices $(g,h), (g_1,h),\ldots, (g_k,h)$ force an induced $Q_{k}$ in the layer $G^h$. Similarly, the vertices $(g,h), (g,h_1),\ldots, (g,h_{d-k})$ force an induced $Q_{d-k}$ in the layer $^gH$. From these two facts and using the unique square property again we conclude that $Q = p_G(Q)\cp p_H(Q)$. \qed

Although the proof of Proposition~\ref{prp:cart-prod-multiplicative} might seems obvious we point out that the conclusion of the proof does not hold for arbitrary subproducts of Cartesian products. Consider for instance the graph $G = P_2\cp P_3$.  Then the product $P_2\cp C_4$ contains subgraphs isomorphic to $G$ that are not products of their projections.

An immediate consequence of Proposition~\ref{prp:cart-prod-multiplicative} is that if $u\in V(Q_n)$, then 
\begin{equation}
\label{eq:distance-cube-poly-for cubes}
D_{Q_n,u}(x,y) = D_{Q_n,0^n}(x,y)=(1+x+y)^n\,,
\end{equation} 
a result earlier obtained in~\cite{saygi-2017+}.

We also point out that the class of daisy cubes is closed under the Cartesian product, a fact which further extends the richness of the class of daisy cubes.

Let $H$ be an induced hypercube of $Q_n$. Then it is well-known that there exists a unique vertex of $H$ with maximum weight, we will call it the {\em top vertex} of $H$ and denote it with $t(H)$. Similarly, $H$ contains a unique vertex with minimum weight to be called the {\em base vertex} of $H$ and denoted $b(H)$. Furthermore $H = \left\langle  I(b(H), t(H))\right\rangle$. 

We are now ready for the main result of this section. 

\begin{theorem}\label{th:DfromC}
If $G$ is a daisy cube, then $D_{G,0^n} (x,y)= C_{G}(x+y-1)$.  
\end{theorem}

\proof
Let $G=Q_n(X)$ and $\widehat{X}=\{x_1,\ldots,x_p\}$ be the maximal vertices of $G$. We thus have $V(G)=\bigcup_{i \in [p]} I(0^n,x_i) $.

An induced $k$-cube $H$ of $Q_n$ is an induced $k$-cube of $G$ if and only if $t(H)\in V(G)$. Similarly an induced $k$-cube $H$ of $Q_n$ is an induced $k$-cube of $ \left\langle I(0^n,x) \right\rangle$ if and only if $t(H)\in I(0^n,x)$.

For any $k$-cube $H$ of $Q_n$ and any subset $T$ of $V(Q_n)$ let $\mathbbm 1_H(T)=1$ if $t(H)\in T$, and $\mathbbm 1_H(T)=0$ otherwise. Let $\mathbb H_{k,d}$ be the set of induced $k$-cubes of $Q_n$ that are at distance $d$ from $0^n$, and let $\mathbb H_k$ be the set of induced $k$-cubes of $Q_n$. Using this notation we have 
\begin{equation}
D_{G,0^n}(x,y)=\sum_k\sum_d\sum_{H \in \mathbb H_{k,d}}\mathbbm 1_H(V(G))x^ky^d
\label{eq:polD}
\end{equation}
and
\begin{equation}
 C_G(z)=\sum_k\sum_{H \in \mathbb H_{k}}\mathbbm 1_H(V(G))z^k\,.
\label{eq:polC}
\end{equation}
Note that if $H$ is a $k$-cube of $Q_n$, then $\mathbbm 1_H(T)=\chi_T(t(H))$, where $\chi_T$ is the characteristic function for the set $T$. By the inclusion-exclusion principle for the union of sets  $A_1,\ldots, A_p$ we thus have 
$$\mathbbm 1_H(\bigcup_{i \in [p]} A_i )=\sum_{J\subset[p],J \neq \emptyset}{(-1)^{|J|-1}\mathbbm 1_H(\bigcap_{i\in J}A_i)}\,.$$
Therefore,
 $$\mathbbm 1_H(V(G)) =\mathbbm 1_H(\bigcup_{i \in [p]} I(0^n,x_i) )=\sum_{J\subset[p],J \neq \emptyset}{(-1)^{|J|-1}\mathbbm 1_H(\bigcap_{i\in J}I(0^n,x_i))}\,.$$ 
Changing the order of summation in~\eqref{eq:polD} and~\eqref{eq:polC} we obtain 
$$ D_{G,0^n}(x,y) =\sum_{J\subset[p],J \neq \emptyset}{(-1)^{|J|-1}D_{\left\langle{\bigcap_{i\in J}I(0^n,x_i)}\right\rangle,0^n}(x,y)}$$ 
and
$$ C_G(z) =\sum_{J\subset[p],J \neq \emptyset}{(-1)^{|J|-1}C_{\left\langle\bigcap_{i\in J}I(0^n,x_i)\right\rangle}(z)}\,.$$
Note that for arbitrary vertices $u,v$ of $Q_n$ we have $I(0^n,u)\cap I(0^n,v)=I(0^n,u\wedge v)$, where $(u\wedge v)_i=1$ if and only if $u_i=1$ and $v_i=1$. The same property extends to the intersection of an arbitrary number of intervals. So $\bigcap_{i\in J}I(0^n,x_i)$ is an interval that induces a hypercube with base vertex $0^n$. From~\eqref{eq:cube-poly-for cubes} and~\eqref{eq:distance-cube-poly-for cubes} we see that the asserted result of the theorem holds if $G$ is an induced hypercube with base vertex $0^n$. Therefore,
$$D_{\left\langle{\bigcap_{i\in J}I(0^n,x_i)}\right\rangle,0^n}(x,y)=C_{\left\langle\bigcap_{i\in J}I(0^n,x_i)\right\rangle}(x+y-1)$$ 
and we are done.
\qed

Theorem~\ref{th:DfromC} has the following immediate consequence.  

\begin{corollary}
\label{cor:commutative} 
If $G$ is a daisy cube, then $D_{G,0^n}(x,y)=D_{G,0^n}(y,x)$.  
\end{corollary}

This corollary in other words says that for any integers $k,d$ the number of induced $k$-cubes at distance $d$ from $0^n$ in a daisy cube $G$ is equal to the number of induced $d$-cubes at distance $k$ from $0^n$.  

To obtain another consequence of Theorem~\ref{th:DfromC} we introduce the counting polynomial of the number of vertices at a given distance from a vertex $u$ as follows. 

\begin{definition}
If $u$ is a vertex of a graph $G$, then let $w_d(G)$, $d\ge 0$, be the number of vertices of $G$ at distance $d$ from $u$. Then set 
$$W_{G,u}(x) = \sum_{d\geq 0} w_{d}(G) x^d\,.$$
\end{definition}

With this definition in hand we can state the following important consequence of Theorem~\ref{th:DfromC}. 

\begin{corollary}
\label{cor:D-and-C-from-W}
If $G$ is a daisy cube, then 
$$D_{G,0^n} (x,y)= W_{G,0^n}(x+y)\quad \mbox{and}\quad C_{G}(x)=W_{G,0^n}(x+1)\,.$$
\end{corollary}

\proof
From Theorem~\ref{th:DfromC} we get $C_G(x) = D_{G,0^n}(0,x+1)$. Since by the definition of the polynomials $D_{G,0^n}(0,x)=W_{G,0^n}(x)$ holds, and consequently $D_{G,0^n}(0,x+1)=W_{G,0^n}(x+1)$, we conclude that $C_G(x) = W_{G,0^n}(x+1)$. 

Using Theorem~\ref{th:DfromC} again and the already proved second assertion of the corollary we get the first assertion: $D_{G,0^n} (x,y)=C_{G}(x+y-1)=W_{G,0^n}(x+y)$. 
\qed

So if $G$ is a daisy cube, then the polynomials $D_{G,0^n}$ and $C_{G}$ are completely determined by $W_{G,0^n}$. 

Consider first hypercubes $Q_n$. Since the number of vertices of weight $k$ in $Q_n$ is $\binom{n}{k}$, we have $W_{Q_n,0^n}(x)=(1+x)^n$. Hence from Corollary~\ref{cor:D-and-C-from-W} we obtain again $C_{Q_n}(x)=(2+x)^n$ and $D_{Q_n,0^n}(x,y)=(1+x+y)^n$.

For the Fibonacci cube $\Gamma_n$ it is well-known that the number of vertices at distance $k$ from $0^n$ is $\binom{n-k+1}{k}$. Therefore $W_{\Gamma_n,0^n}(x)= \sum_{k=0}^{\left\lfloor \frac{n+1}{2}\right\rfloor}{\binom{n-k+1}{k}x^k}$ and we deduce that  
$$C_{\Gamma_n}(x)=\sum_{k=0}^{\left\lfloor \frac{n+1}{2}\right\rfloor}{\binom{n-k+1}{k}(1+x)^k}\,,$$ 
a result first proved in~\cite[Theorem 3.2]{klavzar-2012} and that
\begin{eqnarray*}
D_{\Gamma_n,0^n}(x,y) & = & \sum_{a=0}^{\left\lfloor \frac{n+1}{2}\right\rfloor}{\binom{n-a+1}{a}(x+y)^a} \\
& = & \sum_{k=0}^{\left\lfloor \frac{n+1}{2}\right\rfloor}{\sum_{d=0}^{\left\lfloor \frac{n+1}{2}\right\rfloor}\binom{n-k-d+1}{k+d}\binom{k+d}{d}{x^ky^d}}\,,
\end{eqnarray*}
a result obtained in~\cite[Proposition 3]{saygi-2017+}.

For the Lucas cube $\Lambda_n$ we have $W_{\Lambda_n,0^n}(x)= \sum_{k=0}^{\left\lfloor \frac{n}{2}\right\rfloor}{\left[2\binom{n-k}{k}-\binom{n-k-1}{k}\right]x^k}$. Therefore,   $$C_{\Lambda_n}(x)=\sum_{k=0}^{\left\lfloor \frac{n}{2}\right\rfloor}{\left[2\binom{n-k}{k}-\binom{n-k-1}{k}\right](1+x)^k}\,,$$ 
which is~\cite[Theorem 5.2]{klavzar-2012}. We note in passing that in~\cite[Theorem 5.2]{klavzar-2012} and ~\cite[Corollary 5.3]{klavzar-2012} there is a typo stating $\binom{n-k+1}{k}$ instead of $\binom{n-k-1}{k}$. Moreover, Corollary~\ref{cor:D-and-C-from-W} also gives that
$$D_{\Lambda_n,0^n}(x,y)=\sum_{k=0}^{\left\lfloor \frac{n}{2}\right\rfloor}{\left[2\binom{n-k}{k}-\binom{n-k-1}{k}\right](x+y)^k}\,.$$.

\begin{corollary}
Let $\{G_n\}_{n=0}^{\infty}$ be a family of daisy cubes and let $f(x,z)=\sum_{n\geq 0}W_{G_n,0^n}(x)z^n$ be the generating function of $W_{G_n,0^n}(x)$. Then the generating function of $C_{G_n}(x)$ is  $$g(x,z)=\sum_{n\geq0}C_{G_n}(x)z^n=f(x+1,z)\,,$$ 
and the generating function of $D_{G_n,0^n}(x,y)$ is 
$$h(x,y,z)=\sum_{n\geq0}D_{G_n,0^n}(x,y)z^n=f(x+y,z)\,.$$
\end{corollary}

For example, the respective generating functions for $\{Q_n\}$ are 
$$f(x,z)=\frac{1}{1-z(1+x)},\ g(x,z)=\frac{1}{1-z(2+x)}, h(x,y,z)=\frac{1}{1-z(1+x+y)}\,,$$
and the respective generating functions for $\{\Lambda_n\}$ are 
$$f(x,z)=\frac{1+xz^2}{1-z-xz^2}, g(x,z)=\frac{1+(x+1)z^2}{1-z-(x+1)z^2}, h(x,y,z)=\frac{1+(x+y)z^2}{1-z-(x+y)z^2}\,,$$
cf.~\cite[Proposition 1]{saygi-2017++}. 

\section{A tree-like equality for daisy cubes}
\label{sec:tree-like}

If $G$ is a daisy cube, then the values of $D_{G,u}(x,y)$ and $W_{G,u}(x,y)$ depend of the choice of the vertex $u$. We demonstrate this on the vertex-deleted $3$-cube $Q_3^-$. 

$Q_3^-$ contains three orbits under the action of its automorphism group on the vertex set. Consider their representatives $000$, $100$, and $110$, for which we have the following polynomials:
\begin{itemize}
\item $W_{Q_3^-,000}(x)=1+3x+3x^2$, \\
$D_{Q_3^-,000}(x,y)=1+3y+3y^2+3x+6xy+3x^2$;
\item $W_{Q_3^-,100}(x)=1+3x+2x^2+x^3$, \\
$D_{Q_3^-,100} (x,y) = 1+3y+2y^2+y^3+x(3+4y+2y^2)+x^2(2+y)$; 
\item $W_{Q_3^-,110}(x)=1+2x+3x^2+x^3$, \\
$D_{Q_3^-,110}(x,y) = 1+2y+3y^2+y^3+x(2+4y+3y^2)+x^2(1+2y)$.
\end{itemize}
Note that $D_{Q_3^-,u}(y,x)\neq D_{Q_3^-,u}(x,y)$ except for $u=0^n$. In addition, there is no obvious relation between $D_{Q_3^-,u}$ and $W_{Q_3^-,u}$. On the other hand, $C_{Q_3^-}(x) = 7 + 9x + 3x^2$, and we observe that $D_{Q_3^-,u}(x,-x)=C_{Q_3^-}(-1)=1$ holds for any vertex $u$. 

Recall that a connected graph $G$ is {\em median} if $|I(u,v)\cap I(u,w)\cap I(v,w)|=1$ holds for any triple of vertices $u$, $v$, and $w$ of $G$. It is well-known that median graphs are partial cubes, cf.~\cite[Theorem 2]{klavzar-1999}. Soltan and Chepoi~\cite{soltan-1987} and independently \v{S}krekovski~\cite{skrekovski-2001} proved that if $G$ is a  median graph then $C_G(-1)=1$. This equality in particular generalizes the fact that $n(T) - m(T) = 1$ holds for a tree $T$. Hence if a daisy cube $G$ is median (say a Fibonacci cube), then by Theorem~\ref{th:DfromC} we have $D_{G,0^n}(x,-x)=1$. Our next result (to be proved in the rest of the section) asserts that this equality holds for every daisy cube and every vertex.  

\begin{theorem}
\label{thm:plus1}
If $G$ is a daisy cube, then $D_{G,u}(x,-x)=1 $ holds for every vertex $u$ in $G$.
\end{theorem}

The following consequence extends the class of partial cubes $G$ for which $C_{G}(-1)=1$ holds. 

\begin{corollary}
If $G$ is a daisy cube, then $C_{G}(-1)=1$.
\end{corollary}	

\proof
By Theorem~\ref{th:DfromC} we have $C_{G}(-1)=D_{G,0^n}(x,-x)$. Theorem~\ref{thm:plus1} completes the argument.
\qed

Berrachedi~\cite{berrachedi-1994} characterized median graphs as the connected graphs $G$ in which intervals $I$ have the property that for every $v\in V(G)$, there exists a unique vertex $x\in I$ that attains $d(v,I)$. As $Q_n$ is median and its subgraph $H$ isomorphic to a $k$-cube is an interval of $Q_n$, it follows that if $u\in V(Q_n)$, then $H$ contains a unique vertex that is at the minimum distance from $u$. It will be denoted $\pi_H(u)$. 

Recall that if $u,v\in V(Q_n)$, then $u\wedge v$ is the vertex with $(u\wedge v)_i=1$ if and only if $u_i=1$ and $v_i=1$. The following fact is straightforward. 
    
\begin{lemma}\label{lem:p1}
If $u$ and $v$ are vertices of $Q_n$ and $G=\left\langle I(0^n,v)\right\rangle$, then  $\pi_G(u)=u\wedge v$.
\end{lemma}

\begin{lemma}\label{lem:p2}
Let $c\in V(Q_n)$ and let $u$ be a vertex from $I(0^n,c)$. Then for any vertex $b$ of $Q_n$ we have $\pi_{\left\langle I(0^n,c)\cap I(0^n,b)\right\rangle}(u)=\pi_{\left\langle I(0^n,b)\right\rangle}(u)$.
\end{lemma}

\proof 
Since $u$ belongs to $I(0^n,c)$ we have $u \wedge c=u$ and $I(0^n,c)\cap I(0^n,b)= I(0^n,c \wedge b)$. Therefore, having in mind Lemma~\ref{lem:p1},
\begin{eqnarray*} 
\pi_{\left\langle I(0^n,c)\cap I(0^n,b)\right\rangle}(u) & = & \pi_{\left\langle I(0^n,c \wedge b)\right\rangle}(u) \\
& = & u \wedge c \wedge b = u \wedge b \\
& = & \pi_{\left\langle I(0^n,b)\right\rangle}(u).
\end{eqnarray*}
\qed

Let $u$ be vertex of $Q_n$ and $G$ be a subgraph of $Q_n$. We can consider the distance cube polynomial $D_{G,u} (x,y)$ even if $u$  is not a vertex of $G$. This polynomial is independent of $n$, we need only that $n$ is big enough such that both $u$ is a vertex of $Q_n$ and $G$ is  a subgraph of $Q_n$.

For the next lemma the following concept introduced in~\cite{dress-1987} (see also~\cite{changat-2016}) is useful.  A subgraph $H$ of a graph $G$ is called  {\em gated} if for every $u\in V(G)$, there exists a vertex $x\in V(H)$ such that for every $v\in V(H)$ the vertex $x$ lies on a shortest $u,v$-path. If such a vertex exists, it must be unique.

\begin{lemma}
\label{lem:-x-to-power}
If $u, b\in V(Q_n)$ and $G=\left\langle I(0^n,b)\right\rangle$, then 
$$D_{G,u} (x,-x)= (-x)^{d(u,G)}\,.$$
\end{lemma}

\proof
The contribution of some induced $k$-cube $H$ of $G$  to the polynomial  $D_{G,u}(x,y)$ is $x^ky^\delta$ where $\delta = d(u,H)$. Because intervals are gated in median graphs and hence in hypercubes, see~\cite[Theorem 6(vi)]{klavzar-1999}, we have  
$$\delta= d(u,\pi_G(u))+d(\pi_G(u),H)\,.$$
Therefore $D_{G,u}(x,y)= y^{d(u,G)}D_{G,\pi_G(u)}(x,y)$. Since $G$ is a hypercube, thus a vertex-transitive graph, and $\pi_G(u)$ belongs to $G$, we have $$D_{G,\pi_G(u)}(x,-x)=D_{G,O^n}(x,-x)=1\,.$$  
We conclude that $D_{G,u}(x,-x)=(-x)^{d(u,G)}$.
\qed

\noindent{\bf Proof (of Theorem~\ref{thm:plus1}).}  
Assume that $G=Q_n(\widehat{X})$ with $\widehat{X}=\left\{x_i,i\in I\right\} $ thus $V(G)=\bigcup_{i\in I} I(0^n,x_i)$.
By inclusion-exclusion formula,
 \begin{equation}
\label{ed:ief}
D_{G,u}(x,-x) =\sum_{J \subset I,J \neq \emptyset}{(-1)^{|J|-1}D_{\left\langle\bigcap_{i\in J} I(0^n,x_i)\right\rangle,u}(x,-x)}.
\end{equation}
Since $\bigcap_{i\in J}I(0^n,x_i)$ is some interval $I(0^n,b_J)$, Lemma~\ref{lem:-x-to-power} implies that 
$$D_{\left\langle\bigcap_{i\in J} I(0^n,x_i)\right\rangle,u}(x,-x)=(-x)^{d(u,\left\langle\bigcap_{i\in J}I(0^n,x_i))\right\rangle}\,.$$
Let $i_0 \in I$ such that $u$ belong to $I(0^n,x_{i_0})$ and let $I'=I \setminus \{i_0\}$. For every  $J'$ of $I'$ consider the pair $\{{J',J'\cup\{i_0\}}\} $. We obtain the following partition of the power set $\mathcal{P}(I)$
$$\mathcal{P}(I)=\bigcup_{J'\subset I'}{\{J'\cup(J'\cup\{i_0\})\}}.$$
If $J'$ is not empty then by Lemmas~\ref{lem:p1} and~\ref{lem:p2} we have 
$$d(u,\left\langle\cap_{i\in J'\cup\{i_0\}}I(0^n,x_i)\right\rangle)=d(u,\left\langle\cap_{i\in J'}I(0^n,x_i)\right\rangle)\,,$$ 
and since $|J'\cup\{i_0\}|=|J'|+1$, the sum of contribution in equation (\ref{ed:ief}) of $J'$ and $J'\cup\{i_0\}$ to  $D_{G,u}(x,-x)$ is null. Therefore the only term remaining corresponds to the pair $\{\emptyset,\{i_0\}\}$.
Thus $D_{G,u}(x,-x)=(-x)^{d(u,\left\langle I(0^n,x_{i_0})\right\rangle)}=1$.
\qed



\begin{thebibliography}{50}

\bibitem{albenque-2016}
   M.~Albenque, K.~Knauer, 
   Convexity in partial cubes: the hull number,
   Discrete Math.\ 339 (2016) 866--876.

\bibitem{berrachedi-1994}
  A.~Berrachedi, 
  A new characterization of median graphs, 
  Discrete Math.\ 128 (1994) 385--387. 

\bibitem{bresar-2003}
  B.~Bre\v{s}ar, S.~Klav\v{z}ar, R.~\v{S}krekovski,
  The cube polynomial and its derivatives: the case of median graphs,
  Electron.\ J.\ Combin.\ 10 (2003) \#R3. 

\bibitem{bresar-2006}
  B.~Bre\v{s}ar, S.~Klav\v{z}ar, R.~\v{S}krekovski,
  Roots of cube polynomials of median graphs,
  J.\ Graph Theory 52 (2006) 37--50. 

\bibitem{changat-2016}
  M.~Changat, I.~Peterin, A.~Ramachandran, 
  On gated sets in graphs,
  Taiwanese J.\ Math.\ 20 (2016) 509--521.

\bibitem{deza-2010}
  M.~Deza, M.~Laurent, 
  Geometry of Cuts and Metrics, 
  Springer, Heidelberg, 2010. 

\bibitem{dress-1987}
  A.~W.~M.~Dress, R.~Scharlau, 
  Gated sets in metric spaces, 
  Aequationes Math.\ 34 (1987) 112--120. 

\bibitem{gravier-2015}
  S.~Gravier, M.~Mollard, S.~\v{S}pacapan, S.~S.~Zemlji\v{c}, 
  On disjoint hypercubes in Fibonacci cubes, 
  Discrete Appl.\ Math.\ 190/191 (2015) 50--55.

\bibitem{hammack-2011}
  R.~Hammack, W.~Imrich, S.~Klav\v{z}ar,
  Handbook of Product Graphs. Second Edition,
  CRC Press, Boca Raton, FL, 2011.

\bibitem{ilic-2012}
  A.~Ili\'{c}, S.~Klav\v{z}ar, Y.~Rho,
  Generalized Fibonacci cubes,
  Discrete Math.\ 312 (2012) 2--11.

\bibitem{jukna-2011}
  S.~Jukna, 
  Extremal Combinatorics. With Applications in Computer Science. Second Edition, 
  Springer, Heidelberg, 2011. 

\bibitem{kirlangic-2009}
  A.~Kirlangic, 
  The rupture degree and gear graphs,
  Bull.\ Malays.\ Math.\ Sci.\ Soc.\ 32 (2009) 31--36.

\bibitem{kleitman-1966}
  D.~J.~Kleitman, 
  Families of non-disjoint subsets, 
  J.\ Combin.\ Theory 1 (1966) 153--155. 

\bibitem{klavzar-2013}
  S.~Klav\v{z}ar, 
  Structure of Fibonacci cubes: a survey,
  J.\ Comb.\ Optim.\ 25 (2013) 505--522.

\bibitem{klavzar-2012}
  S.~Klav\v{z}ar, M.~Mollard,
  Cube polynomial of Fibonacci and Lucas cubes,
  Acta Appl.\ Math.\ 117 (2012) 93--105. 

\bibitem{klavzar-1999}
  S.~Klav\v{z}ar, H.~M.~Mulder, 
  Median graphs: characterizations, location theory and related structures,
  J.\ Combin.\ Math.\ Combin.\ Comput.\ 30 (1999) 103--127.  

\bibitem{liu-1994}
  J.~Liu, W.-J.~Hsu, M.~J.~Chung,
  Generalized Fibonacci cubes are mostly Hamiltonian,
  J.\ Graph Theory 18 (1994) 817--829.

\bibitem{marc-2016}
  T.~Marc, 
  There are no finite partial cubes of girth more than 6 and minimum degree at least 3,
  European J.\ Combin.\ 55 (2016) 62--72.

\bibitem{mollard-2017}
  M.~Mollard, 
  Non covered vertices in Fibonacci cubes by a maximum set of disjoint hypercubes,
  Discrete Appl.\ Math.\ 219 (2017) 219--221.

\bibitem{munarini-2001}
  E.~Munarini, C.~Perelli Cippo, N.~Zagaglia Salvi, 
  On the Lucas cubes,
  Fibonacci Quart.\ 39 (2001) 12--21.

\bibitem{ovchinnikov-2011}
  S.~Ovchinnikov, 
  Graphs and Cubes,
  Springer, New York, 2011. 

\bibitem{saygi-2016}
  E.~Sayg{\i}, \"O.~E\u{g}ecio\u{g}lu, 
  Counting disjoint hypercubes in Fibonacci cubes, 
  Discrete Appl.\ Math.\ 215 (2016) 231--237.

\bibitem{saygi-2017+}
  E.~Sayg{\i}, \"O.~E\u{g}ecio\u{g}lu, 
  $q$-cube enumerator polynomial of Fibonacci cubes, 
  Discrete Appl.\ Math.\ (2017), in press, doi.org/10.1016/j.dam.2017.04.026.  

\bibitem{saygi-2017++}
  E.~Sayg{\i}, \"O.~E\u{g}ecio\u{g}lu, 
  $q$-counting hypercubes in Lucas cubes,
  Turkish J.\ Math., in press.  

\bibitem{skrekovski-2001}
  R.~\v{S}krekovski, 
  Two relations for median graphs, 
  Discrete Math.~ 226 (2001) 351--353.

\bibitem{soltan-1987}
  P.~S.~Soltan, V.~D.~Chepoi, 
  Solution of the Weber problem for discrete median metric spaces, 
  Trudy Tbiliss.\ Mat.\ Inst.\ 85 (1987) 53-76 (in Russian).

\bibitem{taranenko-2013}
  A.~Taranenko, 
  A new characterization and a recognition algorithm of Lucas cubes,
  Discrete Math.\ Theor.\ Comput.\ Sci.\ 15 (2013) 31--39. 

\bibitem{wei-2016}
  J.~Wei, H.~Zhang, 
  Proofs of two conjectures on generalized Fibonacci cubes,
  European J.\ Combin.\ 51 (2016) 419--432.

\end{thebibliography}
\end{document}